\documentclass[11pt]{article}
\usepackage{amssymb}
\usepackage{amsfonts}
\usepackage{amsmath}

\input{tcilatex}
\begin{document}

\begin{center}
\textbf{ON THE LIPSCHITZ STABILITY OF INVERSE NODAL PROBLEM FOR DIRAC SYSTEM}

\textbf{Emrah YILMAZ\ and Hikmet KOYUNBAKAN }

\textbf{Firat University, Department of Mathematics, 23119, Elaz\i g / TURKEY%
}

\textbf{emrah231983@gmail.com\qquad hkoyunbakan@gmail.com}
\end{center}

\begin{quote}
\textbf{Abstract: }{\footnotesize Inverse nodal problem on Dirac operator is
finding the parameters\ in the boundary conditions, number }$m$ 
{\footnotesize and potential function }$V${\footnotesize \ by using a set of
nodal points of a component\ of two component vector eigenfunctions as the
given spectral data. In this study, we solve a stability problem using nodal
set of vector eigenfunctions and show that the space of all }$V$%
{\footnotesize \ functions is homeomorphic to the partition set of all space
of asymptotically equivalent nodal sequences induced by an equivalence
relation. Moreover, we give a reconstruction formula for the potential
function as a limit of a sequence of functions and associated nodal data of
one component of vector eigenfunction. Our technique depends on the explicit
asymptotic expressions of the nodal parameters and, it is basically similar
to} {\footnotesize \cite{Tsay,Cheng} which is given for Sturm-Liouville and
Hill's operators, respectively.}

\smallskip



\smallskip
\end{quote}

\textbf{\large 1. Introduction}

Inverse spectral problems have been a significant research area in
mathematical physics. Different methods have been proposed to recover
coefficient functions in differential equations by using spectral data \cite%
{AMB,Lev,Lagh,Pschl,PVVRCHK,BT,Run,Pronska}$.$ Generally, the spectral data
have consisted of the eigenvalues and a corresponding sequence of norming
constants, or two eigenvalue sequences. In 1988, McLaughlin showed that
knowledge of nodal points can determine the potential function of
Sturm-Liouville problem up to a constant \cite{Mc}$.$ This is so called
inverse nodal problem. Numerical schemes were then given by Hald and
McLaughlin \cite{Hald} to reconstruct the density function of a vibrating
string, the elastic modulus of a vibrating rod, the potential function in
Sturm-Liouville problem. Independently, Shen \textit{et al.} \cite{Shen}
studied the relation between nodal points and density function of string
equation in 1988. Many results and reconstruction formulas have been derived
about inverse nodal problem by several authors \cite%
{sleeman,Kryshva,Koyunbakan,guo,Law1}. Here, we deal with the inverse nodal
problem for Dirac system.

Dirac system is a modern presentation of the relativistic quantum mechanics
of electrons intended to make new mathematical results accesible to a wider
audience. It treats in some depth relativistic invariance of a quantum
theory, self-adjointness and spectral theory, qualitative features of
relativistic bound and scattering states and the external field problem in
quantum electrodynamics, without neglecting the interpretational
difficulties and limitations of the theory \cite{emrah}.

Inverse problems for Dirac system had been investigated by Moses \cite{Moses}%
, Prats and Toll \cite{prats}, Verde \cite{vERDE}, Gasymov and Levitan \cite%
{GASY}, and Panakhov \cite{PAN}. It is well known that two spectra uniquely
determine the matrix valued potential function in Dirac system \cite{DZH}.
In \cite{JOA}, eigenfunction expansions for one dimensional Dirac operator
describing the motion of a particle in quantum mechanics were investigated.
In addition, inverse spectral problems for weighted Dirac system were
studied in \cite{WATSON}.

One studied the properties of the eigenvalues and vector-valued
eigenfunctions for the Dirac system with the same spectral parameter in the
equations and the boundary conditions \cite{Kerim}. Sampling theory of
signal analysis associated with Dirac systems, when the eigenvalue parameter
appears linearly in the boundary conditions was investigated in \cite{Annaby}%
. One investigated a problem for the Dirac differential operators in the
case where an eigenparameter not only appears in the differential equation
but is also linearly contained in a boundary condition, and proved
uniqueness theorems for inverse spectral problem with known collection of
eigenvalues and normalizing constants or two spectra \cite{Amirov}. Other
than these studies, there are many papers in literature (see \cite%
{MYKY1,MYKY2,MYKY3,Thaller,Puyda}).

Inverse nodal problems for Dirac system had not been studied until the works
of Yang and Huang \cite{YANGHUAN}. They gave reconstruction formulas for one
dimensional Dirac operator by using nodal datas. Later years, inverse nodal
problem was solved for Dirac system under different boundary conditions \cite%
{Tuba,YANGPV}.

Consider the Dirac system 
\begin{equation}
By^{\prime }(x)+Q(x)y(x)=\lambda y(x),0\leq x\leq \pi ,  \tag{1.1}
\end{equation}%
with boundary conditions%
\begin{eqnarray}
\left( \lambda \cos \alpha +a_{0}\right) y_{1}(0)+\left( \lambda \sin \alpha
+b_{0}\right) y_{2}(0) &=&0,  \notag \\
\left( \lambda \cos \beta +a_{1}\right) y_{1}(\pi )+\left( \lambda \sin
\beta +b_{1}\right) y_{2}(\pi ) &=&0,  \TCItag{1.2}
\end{eqnarray}%
where $\lambda $ is a spectral parameter,%
\begin{equation}
\begin{array}{ccc}
B=\left( 
\begin{array}{cc}
0 & 1 \\ 
-1 & 0%
\end{array}%
\right) , & Q(x)=\left( 
\begin{array}{cc}
V(x)+m & 0 \\ 
0 & V(x)-m%
\end{array}%
\right) , & y(x)=\left( 
\begin{array}{c}
y_{1}(x) \\ 
y_{2}(x)%
\end{array}%
\right) ,%
\end{array}
\tag{1.3}
\end{equation}%
and $V$ is a real valued, continuous function on $[0,\pi ].$ Furthermore, $%
m,a_{k},b_{k}(k=0,1),\alpha $ and $\beta $ are real constants: moreover $-%
\frac{\pi }{2}\leq \alpha ,\beta \leq \frac{\pi }{2}$ \cite{YANGPV}$.$
Throughout the paper \cite{YANGPV}, Yang and Pivovarchik supposed that%
\begin{eqnarray}
a_{0}\sin \alpha -b_{0}\cos \alpha &>&0\text{,}  \notag \\
a_{1}\sin \beta -b_{1}\cos \beta &<&0.  \TCItag{1.4}
\end{eqnarray}%
The properties of the eigenvalues and eigenfunctions of the problem
(1.1)-(1.2) were studied in \cite{Kerim}. Under the condition (1.4), the
eigenvalues of the problem (1.1)-(1.2) are real and algebraically simple 
\cite{Kerim}. Considering (1.3) in (1.1), we get 
\begin{equation*}
\left( 
\begin{array}{cc}
0 & 1 \\ 
-1 & 0%
\end{array}%
\right) \left( 
\begin{array}{c}
y_{1}^{\prime }(x) \\ 
y_{2}^{\prime }(x)%
\end{array}%
\right) +\left( 
\begin{array}{cc}
V(x)+m & 0 \\ 
0 & V(x)-m%
\end{array}%
\right) \left( 
\begin{array}{c}
y_{1}(x) \\ 
y_{2}(x)%
\end{array}%
\right) =\lambda \left( 
\begin{array}{c}
y_{1}(x) \\ 
y_{2}(x)%
\end{array}%
\right) ,
\end{equation*}%
and thus, equation (1.1) is equivalent to a system of two simultaneous first
order differential equations%
\begin{eqnarray}
y_{2}^{\prime }(x,\lambda )+\left[ V(x)+m\right] y_{1}(x,\lambda )-\lambda
y_{1}(x,\lambda ) &=&0,  \notag \\
-y_{1}^{\prime }(x,\lambda )+\left[ V(x)-m\right] y_{2}(x,\lambda )-\lambda
y_{2}(x,\lambda ) &=&0.  \TCItag{1.5}
\end{eqnarray}%
In general, potential function of Dirac system (1.1) has the following form%
\begin{equation*}
Q(x)=\left( 
\begin{array}{cc}
p_{11}(x) & p_{12}(x) \\ 
p_{21}(x) & p_{22}(x)%
\end{array}%
\right) ,
\end{equation*}%
where $p_{ik}(x)$ $(i,k=1,2)$ are real valued and continuous functions on $%
[0,\pi ].$ For the case in which $p_{12}(x)=p_{21}(x)=0$ and $%
p_{11}(x)=V(x)+m,$ $p_{22}(x)=V(x)-m$ where $m$ is the mass of particle, the
system (1.5) is known in relativistic quantum theory as a stationary one
dimensional Dirac system or first canonical form of Dirac system \cite{Lev}.

Let $y(x,\lambda _{n})=\left[ y_{1}(x,\lambda _{n}),y_{2}(x,\lambda _{n})%
\right] ^{T}$ be two dimensional vector eigenfunction of the Dirac system
(1.1) related to the eigenvalue $\lambda =\lambda _{n}$ where $T$ denotes
transpose. Assume that $x_{n}^{j,i}$ are the nodal points of $i-$th
component $y_{i}(x,\lambda _{n})$ of the $n-$th eigenfunction $y(x,\lambda
_{n})$ where $0<x_{n}^{1,i}<x_{n}^{2,i}<...<x_{n}^{n-1,i}<\pi .$ In other
words, $y_{i}(x_{n}^{j,i},\lambda _{n})=0.$ Let $I_{n}^{j,i}=\left(
x_{n}^{j,i},x_{n}^{j+1,i}\right) $ be the $j-$th nodal domain$,$ and let%
\begin{equation*}
l_{n}^{j,i}=x_{n}^{j+1,i}-x_{n}^{j,i},
\end{equation*}%
be the associated nodal length. For simplicity, we agree that $x_{n}^{0,i}=0$
and $x_{n}^{\left\vert n\right\vert +1-i,i}=\pi .$ We also define the
function $j_{n,i}(x)$ to be the largest index $j_{i}$ such that $0\leq
x_{n}^{j,i}\leq x$ for $n>0$ and $j_{n,i}(x)$ to be the largest index $j_{i}$
such that $0\leq x\leq x_{n}^{j,i}$ for $n<0.$ Thus, $j_{i}=j_{n,i}(x)$ if
and only if $x\in \lbrack x_{n}^{j,i},x_{n}^{j+1,i})$ for $n>0$ and $%
(x_{n}^{j+1,i},x_{n}^{j,i}]$ for $n<0$ \cite{YANGHUAN}$.$

Denote $\Lambda ^{i}=\{x_{n}^{j,i}\},i=1,2.$ Hence, $\Lambda =\Lambda
^{1}\cup \Lambda ^{2}$ is called the set of all nodal points of Dirac
operator. This set is dense on $[0,\pi ]$ \cite{YANGPV}$.$ Throughout this
study, we'll give all proofs for the first component of the eigenfunction.

The rest of this study is arranged as follows: in remaining part of section
1, we give some properties of Dirac system and quote some important results
to use in main theorems. In section 2, we obtain some reconstruction
formulas for potential function under different boundary conditions.
Finally, we define $d_{0}$, $d_{\Sigma _{Dir}}$ to prove Lipschitz stability
of inverse nodal problem. Then, we express Theorem 3.1 in section 3.

Now, we need to remind some conclusions which are given by \cite{YANGPV} to
use in our main results.

\bigskip

{\large \textbf{Lemma 1.1.}} \cite{YANGPV} \textit{The spectrum of the
problem (1.1)-(1.2) consists of eigenvalues }$\left\{ \lambda _{n}\right\}
_{n\in 
\mathbb{Z}
\text{ }}$ \textit{which are all real and algebraically simple behave
asymptotically as}%
\begin{equation}
\lambda _{n}=n-2+\frac{v}{\pi }+\frac{c}{n}+O\left( \frac{1}{n^{2}}\right)
,n\rightarrow \infty ,  \tag{1.6}
\end{equation}%
\textit{and}%
\begin{equation*}
\lambda _{-n}=-n+\frac{v}{\pi }-\frac{c}{n}+O\left( \frac{1}{n^{2}}\right)
,n\rightarrow \infty ,
\end{equation*}%
\textit{where}%
\begin{equation*}
v=\dint\limits_{0}^{\pi }V(t)dt+\beta -\alpha ,c=\frac{m^{2}}{2}+\frac{m}{%
2\pi }\left( \sin 2\alpha -\sin 2\beta \right) +\frac{a_{0}\sin \alpha
-b_{0}\cos \alpha }{\pi }-\frac{a_{1}\sin \beta -b_{1}\cos \beta }{\pi }.
\end{equation*}

\bigskip

{\large \textbf{Lemma 1.2. }}\cite{YANGPV} \textit{Let }$y(x,\lambda
)=\left( 
\begin{array}{c}
y_{1}(x,\lambda ) \\ 
y_{2}(x,\lambda )%
\end{array}%
\right) $ \textit{be\ the solution of (1.1) satisfying the condition}%
\begin{equation}
y(0,\lambda )=\left( 
\begin{array}{c}
-\left( \lambda \sin \alpha +b_{0}\right) \\ 
\lambda \cos \alpha +a_{0}%
\end{array}%
\right) ,  \tag{1.7}
\end{equation}%
\textit{then, we have} 
\begin{eqnarray}
y_{1}(x,\lambda ) &=&-\lambda \sin \left( \lambda
x-\dint\limits_{0}^{x}V(t)dt+\alpha \right) +\frac{m^{2}}{2}x\cos \left(
\lambda x-\dint\limits_{0}^{x}V(t)dt+\alpha \right)  \notag \\
&&-m\cos \alpha \sin \left( \lambda x-\dint\limits_{0}^{x}V(t)dt\right)
-a_{0}\sin \left( \lambda x-\dint\limits_{0}^{x}V(t)dt\right)  \TCItag{1.8}
\\
&&-b_{0}\cos \left( \lambda x-\dint\limits_{0}^{x}V(t)dt\right) +O\left( 
\frac{e^{\tau x}}{\lambda }\right) ,  \notag
\end{eqnarray}%
\textit{and}%
\begin{eqnarray}
y_{2}(x,\lambda ) &=&\lambda \cos \left( \lambda
x-\dint\limits_{0}^{x}V(t)dt+\alpha \right) +\frac{m^{2}}{2}x\sin \left(
\lambda x-\dint\limits_{0}^{x}V(t)dt+\alpha \right)  \notag \\
&&+m\sin \alpha \sin \left( \lambda x-\dint\limits_{0}^{x}V(t)dt\right)
+a_{0}\cos \left( \lambda x-\dint\limits_{0}^{x}V(t)dt\right)  \TCItag{1.9}
\\
&&-b_{0}\sin \left( \lambda x-\dint\limits_{0}^{x}V(t)dt\right) +O\left( 
\frac{e^{\tau x}}{\lambda }\right) ,  \notag
\end{eqnarray}%
\textit{where} $\tau =\left\vert \func{Im}\lambda \right\vert .$

\bigskip

{\large \textbf{Lemma 1.3. }} \cite{YANGPV} \textit{For sufficiently large} $%
n>0,$ \textit{the first component} $y_{1}(x,\lambda _{n})$ \textit{of the
eigenfunction} $y(x,\lambda _{n})$ \textit{for Dirac system has exactly }$%
N(\alpha ,\beta )$ \textit{nodes in the interval} $(0,\pi )$ \textit{where}%
\begin{equation*}
N(\alpha ,\beta )=\left\{ 
\begin{array}{ccc}
n-2, & \text{\textit{for} }\alpha \geq 0\text{ \textit{and} }\beta >0\text{ 
\textit{or}} & \text{\textit{for} }\alpha <0\text{ \textit{and} }\beta \leq 0
\\ 
n-3, & \text{\textit{for} }\alpha \geq 0\text{ \textit{and} }\beta \leq 0 & 
\\ 
n-1, & \text{\textit{for} }\alpha <0\text{ \textit{and} }\beta >0. & 
\end{array}%
\right.
\end{equation*}%
Moreover, uniformly with respect to $j\in \{1,2,...,N(\alpha ,\beta )\}$,
the nodal parameters of the problem (1.1)-(1.7) has the following asymptotic
formulas, respectively for sufficiently large $n,$%
\begin{equation*}
x_{n}^{j,1}=\frac{2\lambda _{n}^{2}}{2\lambda _{n}^{2}-m^{2}}\left[ \frac{%
j\pi }{\lambda _{n}}+\frac{1}{\lambda _{n}}\dint%
\limits_{0}^{x_{n}^{j,1}}V(t)dt-\frac{\alpha }{\lambda _{n}}+\frac{m\sin
2\alpha }{2\lambda _{n}^{2}}+\frac{a_{0}\sin \alpha -b_{0}\cos \alpha }{%
\lambda _{n}^{2}}+O\left( \frac{1}{\lambda _{n}^{3}}\right) \right] ,
\end{equation*}%
and%
\begin{equation*}
l_{n}^{j,1}=\frac{2\lambda _{n}^{2}}{2\lambda _{n}^{2}-m^{2}}\left[ \frac{%
\pi }{\lambda _{n}}+\frac{1}{\lambda _{n}}\dint%
\limits_{x_{n}^{j,1}}^{x_{n}^{j+1,1}}V(t)dt+O\left( \frac{1}{\lambda _{n}^{3}%
}\right) \right] ,
\end{equation*}%
where $n\neq \mp \frac{m}{\sqrt{2}}+2.$ Now, we consider the system (1.1)
with boundary conditions%
\begin{eqnarray}
u_{1}(0)\cos \widetilde{\alpha }+u_{2}(0)\sin \widetilde{\alpha } &=&0, 
\notag \\
u_{1}(\pi )\cos \widetilde{\beta }+u_{2}(\pi )\sin \widetilde{\beta } &=&0, 
\TCItag{1.10}
\end{eqnarray}%
where $0\leq \widetilde{\alpha },\widetilde{\beta }\leq \pi ,$ and $m$ is
positive in (1.3). It is well known that the spectrum of the system (1.1)
with the boundary conditions (1.10) includes the eigenvalues $\widetilde{%
\lambda }_{n},n\in 
\mathbb{Z}
$ which are all real and simple, and the sequence $\{\widetilde{\lambda }%
_{n}\}$ satisfies the classical asymptotic form \cite{Lev}, \cite{YANGHUAN}%
\begin{equation}
\widetilde{\lambda }_{n}=n+\frac{\widetilde{v}}{\pi }+\frac{\widetilde{c}_{1}%
}{n}+O\left( \frac{1}{n^{2}}\right) ,  \tag{1.11}
\end{equation}%
where%
\begin{equation*}
\widetilde{v}=\widetilde{\beta }-\widetilde{\alpha }+\dint\limits_{0}^{\pi }%
\widetilde{V}(t)dt,\widetilde{c}_{1}=\frac{m(\sin 2\widetilde{\alpha }-\sin 2%
\widetilde{\beta })+m^{2}\pi }{2\pi \cos ^{2}\left( \dint\limits_{0}^{\pi }%
\widetilde{V}(t)dt-\widetilde{\alpha }+\widetilde{\beta }\right) }.
\end{equation*}

Let $u(x,\widetilde{\lambda })=\left( u_{1}(x,\widetilde{\lambda }),u_{2}(x,%
\widetilde{\lambda })\right) ^{T}$ be the solution of the system (1.1) with
initial conditions%
\begin{equation}
u_{1}(0,\widetilde{\lambda })=\sin \widetilde{\alpha },u_{2}(0,\widetilde{%
\lambda })=-\cos \widetilde{\alpha }.  \tag{1.12}
\end{equation}%
Then, by successive approximations method, there hold%
\begin{eqnarray}
u_{1}(x,\widetilde{\lambda }) &=&\sin \left( \widetilde{\lambda }%
x-\dint\limits_{0}^{x}\widetilde{V}(t)dt+\widetilde{\alpha }\right) -\frac{%
U_{1}}{\widetilde{\lambda }}+O\left( \frac{e^{\left\vert \tau \widetilde{%
\lambda }\right\vert x}}{\widetilde{\lambda }^{2}}\right) ,  \notag \\
u_{2}(x,\widetilde{\lambda }) &=&-\cos \left( \widetilde{\lambda }%
x-\dint\limits_{0}^{x}\widetilde{V}(t)dt+\widetilde{\alpha }\right) -\frac{%
U_{2}}{\widetilde{\lambda }}+O\left( \frac{e^{\left\vert \tau \widetilde{%
\lambda }\right\vert x}}{\widetilde{\lambda }^{2}}\right) ,  \TCItag{1.13}
\end{eqnarray}%
for large $\left\vert \widetilde{\lambda }\right\vert ,$ where \cite%
{YANGHUAN}%
\begin{eqnarray}
U_{1}(x,\widetilde{\lambda }) &=&-m\sin \left( \widetilde{\lambda }%
x-\dint\limits_{0}^{x}\widetilde{V}(t)dt\right) \cos \widetilde{\alpha }+%
\frac{m^{2}}{2}x\cos \left( \widetilde{\lambda }x-\dint\limits_{0}^{x}%
\widetilde{V}(t)dt+\widetilde{\alpha }\right) ,  \notag \\
U_{2}(x,\widetilde{\lambda }) &=&m\sin \left( \widetilde{\lambda }%
x-\dint\limits_{0}^{x}\widetilde{V}(t)dt\right) \sin \widetilde{\alpha }+%
\frac{m^{2}}{2}x\sin \left( \widetilde{\lambda }x-\dint\limits_{0}^{x}%
\widetilde{V}(t)dt+\widetilde{\alpha }\right) .  \TCItag{1.14}
\end{eqnarray}

\bigskip

{\large \textbf{Lemma 1.4. }} \cite{YANGHUAN} \textit{For sufficiently large 
}$\left\vert n\right\vert ,$ \textit{the} $i-$\textit{th component} $u_{i}(x,%
\widetilde{\lambda }_{n})$ \textit{of the eigenfunction} $u(x,\widetilde{%
\lambda }_{n})$ \textit{of the problem} \textit{(1.1),(1.12) has exactly} $%
\left\vert n\right\vert +1-i$ \textit{nodes in the interval} $(0,\pi ).$ 
\textit{Moreover, the asymptotic formulas for nodal points of first and
second components of the eigenfunction} $u(x,\widetilde{\lambda }_{n})$ 
\textit{as} $\left\vert n\right\vert \rightarrow \infty $ \textit{uniformly
with respect to} $j\in 
\mathbb{Z}
$ \textit{are as following}%
\begin{equation}
\widetilde{x}_{n}^{j,1}=\frac{2\widetilde{\lambda }_{n}^{2}}{2\widetilde{%
\lambda }_{n}^{2}-(-1)^{j}m^{2}}\left[ \frac{j\pi }{\widetilde{\lambda }_{n}}%
+\frac{1}{\widetilde{\lambda }_{n}}\dint\limits_{0}^{\widetilde{x}_{n}^{j,1}}%
\widetilde{V}(t)dt-\frac{\widetilde{\alpha }}{\widetilde{\lambda }_{n}}+%
\frac{(-1)^{j}m\sin 2\widetilde{\alpha }}{2\widetilde{\lambda }_{n}^{2}}%
+O\left( \frac{1}{\widetilde{\lambda }_{n}^{3}}\right) \right] ,  \tag{1.15}
\end{equation}%
\textit{and}%
\begin{equation}
\widetilde{x}_{n}^{j,2}=\frac{2\widetilde{\lambda }_{n}^{2}}{2\widetilde{%
\lambda }_{n}^{2}-(-1)^{j}m^{2}}\left[ \frac{\left( j-\frac{1}{2}\right) \pi 
}{\widetilde{\lambda }_{n}}+\frac{1}{\widetilde{\lambda }_{n}}%
\dint\limits_{0}^{\widetilde{x}_{n}^{j,2}}\widetilde{V}(t)dt-\frac{%
\widetilde{\alpha }}{\widetilde{\lambda }_{n}}+\frac{(-1)^{j+1}m\sin 2%
\widetilde{\alpha }}{2\widetilde{\lambda }_{n}^{2}}+O\left( \frac{1}{%
\widetilde{\lambda }_{n}^{3}}\right) \right] ,  \tag{1.16}
\end{equation}%
\textit{where} $n\neq \mp \frac{(-1)^{\frac{j}{2}}m}{\sqrt{2}}$. \textit{The
nodal lengths} $\widetilde{l}_{n}^{j,1}$ \textit{for the problem
(1.1),(1.12) have the following asymptotic expansions}%
\begin{equation*}
\widetilde{l}_{n}^{j,1}=\frac{\pi }{\widetilde{\lambda }_{n}}+\frac{1}{%
\widetilde{\lambda }_{n}}\dint\limits_{\widetilde{x}_{n}^{j,1}}^{\widetilde{x%
}_{n}^{j+1,1}}\widetilde{V}(t)dt+\frac{\{(-1)^{j+1}-(-1)^{j}\}m\sin 2%
\widetilde{\alpha }}{2\widetilde{\lambda }_{n}^{2}}+\frac{\{(-1)^{j+1}%
\widetilde{x}_{n}^{j+1,1}-(-1)^{j}\widetilde{x}_{n}^{j,1}\}m^{2}}{2%
\widetilde{\lambda }_{n}^{2}}+O\left( \frac{1}{\widetilde{\lambda }_{n}^{3}}%
\right) .
\end{equation*}

\textit{In case of} $j=2k$ $($\textit{or} $j=2k+1),k\in 
\mathbb{Z}
;$ \textit{we get}%
\begin{equation}
\widetilde{l}_{n}^{j,1}=\frac{2\widetilde{\lambda }_{n}^{2}}{2\widetilde{%
\lambda }_{n}^{2}\pm m^{2}}\left[ \frac{\pi }{\widetilde{\lambda }_{n}}+%
\frac{1}{\widetilde{\lambda }_{n}}\dint\limits_{\widetilde{x}_{n}^{j,1}}^{%
\widetilde{x}_{n}^{j+1,1}}\widetilde{V}(t)dt\pm \frac{m\sin 2\widetilde{%
\alpha }}{2\widetilde{\lambda }_{n}^{2}}+O\left( \frac{1}{\widetilde{\lambda 
}_{n}^{3}}\right) \right] .  \tag{1.17}
\end{equation}

We can easily obtain $\widetilde{l}_{n}^{j,2}$ similarly as $\left\vert
n\right\vert \rightarrow \infty $ by using definition of nodal lengths and
(1.16). Here, $\{x_{n}^{j,i}\}$, $\{\widetilde{x}_{n}^{j,i}\},i=1,2$ and $%
\{\lambda _{n}\},\{\widetilde{\lambda }_{n}\}$ are the nodal sets and
eigenvalues of the problems (1.1), (1.7) and (1.1), (1.12), respectively.

\bigskip

{\large \textbf{Theorem 1.1.}} \textit{Suppose that }$V\in L_{1}(0,\pi ).$ 
\textit{Then, for almost every} $x\in (0,\pi ),$ \textit{with} $%
j_{i}=j_{n,i}(x),$%
\begin{equation*}
\begin{array}{cc}
\lim\limits_{n\rightarrow \infty }\lambda
_{n}\dint\limits_{x_{n}^{j,i}}^{x_{n}^{j+1,i}}V(t)dt=V(x), & 
\lim\limits_{n\rightarrow \infty }\lambda
_{n}\dint\limits_{x_{n}^{j,i}}^{x_{n}^{j+1,i}}\cos (2\lambda _{n}\pi
t)V(t)dt=0,%
\end{array}%
\end{equation*}%
\textit{where} $i=1,2$ \textit{and }$\lambda _{n}=n-2$ \textit{for the
problem (1.1), (1.7). We can express the similar theorem for the problem
(1.1), (1.12).}

\bigskip

{\large \textbf{Proof: }}It can be proved by similar method given in \cite%
{Tsay}.

\bigskip

{\large \textbf{Remark: }} \cite{YANGHUAN} $\{x_{n}^{j,i}\}\subset \Lambda
^{i},$ $\{\widetilde{x}_{n}^{j,i}\}\subset \widetilde{\Lambda }^{i}$ are
chosen such that 
\begin{equation*}
\lim_{n\rightarrow \infty }x_{n}^{j,i}=x=\lim_{n\rightarrow \infty }%
\widetilde{x}_{n}^{j,i},
\end{equation*}%
where $i=1,2$ and $x\in \lbrack 0,\pi ].$

\bigskip

\textbf{\large 2. Reconstruction of potential function by using nodal points 
}

\bigskip

In this section, we will derive some reconstruction formulas of potential
functions $V$ and $\widetilde{V}$ where $l_{n}^{j,i},\widetilde{l}_{n}^{j,i}$
and $x_{n}^{j,i},\widetilde{x}_{n}^{j,i}$ are nodal lengths and nodal points
for the problems (1.1), (1.7) and (1.1), (1.12), respectively. Here, all of
our proofs and definitions will be given for the first component of
eigenfunction $($That is, for $i=1)$.

\bigskip

{\large \textbf{Theorem 2.1.}} \textit{Let }$V,\widetilde{V}\in L_{1}[0,\pi
] $\textit{\ be the potential functions for Dirac system under the
conditions (1.7) and (1.12), respectively. Define }$F_{n}$\textit{\ by}

\textit{a) For the problem (1.1), (1.7),}%
\begin{equation*}
F_{n}(x)=(n-2)\left\{ \dsum\limits_{j=1}^{n-1}\left[ n-2-\frac{m^{2}}{2(n-2)}%
\right] l_{n}^{j,1}-\pi \right\} .
\end{equation*}

\textit{b) For the problem (1.1), (1.12),}%
\begin{equation*}
\widetilde{F}_{n}(x)=n\left\{ \dsum\limits_{j=1}^{n-1}\left[ n\pm \frac{m^{2}%
}{2n}\right] \widetilde{l}_{n}^{j,i}\pm \frac{m^{2}}{n}\widetilde{x}%
_{n}^{j,i}-\pi \right\} .
\end{equation*}%
\textit{where} $l_{n}^{j,i}=x_{n}^{j+1,i}-x_{n}^{j,i}$ and $\widetilde{l}%
_{n}^{j,i}=\widetilde{x}_{n}^{j+1,i}-\widetilde{x}_{n}^{j,i}$ $(j$ \textit{%
is odd or even})$.$ \textit{Then}, $F_{n}$ and $\widetilde{F}_{n}$ \textit{%
converge to} $V$ and $\widetilde{V}$ \textit{pointwisely almost everywhere,
respectively and also in} $L_{1}$ \textit{sense}.\textit{\ Moreover,
pointwise convergence holds for all the continuity points of} $V$ \textit{and%
} $\widetilde{V}.$

{\large \textbf{Proof:}}

\textit{a) }We will consider the reconstruction formula for the potential
function of the problem (1.1), (1.7). Observe that, by Lemma 1.3, we have 
\begin{equation*}
\frac{\lambda _{n}}{\pi }l_{n}^{j,1}-1-\frac{m^{2}l_{n}^{j,1}}{2\lambda
_{n}\pi }=\frac{1}{\pi }\dint\limits_{x_{n}^{j,1}}^{x_{n}^{j+1,1}}V(t)dt+O%
\left( \frac{1}{\lambda _{n}^{2}}\right) ,
\end{equation*}%
and%
\begin{equation*}
\lambda _{n}\left[ l_{n}^{j,1}\left( \lambda _{n}-\frac{m^{2}}{2\lambda _{n}}%
\right) -\pi \right] =\lambda
_{n}\dint\limits_{x_{n}^{j,1}}^{x_{n}^{j+1,1}}V(t)dt+O\left( \frac{1}{%
\lambda _{n}}\right) .
\end{equation*}

Then, by using asymptotic expansions for eigenvalues, we obtain%
\begin{eqnarray*}
\lambda _{n}\left[ l_{n}^{j,1}\left( \lambda _{n}-\frac{m^{2}}{2\lambda _{n}}%
\right) -\pi \right] &=&\left[ n-2+O\left( \frac{1}{n}\right) \right] \left[
l_{n}^{j,1}\left( n-2+O\left( \frac{1}{n}\right) -\frac{m^{2}}{%
2(n-2)+O\left( \frac{1}{n}\right) }\right) -\pi \right] \\
&=&\left[ n-2+O\left( \frac{1}{n}\right) \right] \left[ l_{n}^{j,1}\left(
n-2-\frac{m^{2}}{2(n-2)}+O\left( \frac{1}{n}\right) \right) -\pi \right] \\
&=&(n-2)\left[ l_{n}^{j,1}\left( n-2-\frac{m^{2}}{2(n-2)}\right) -\pi \right]
+o(1).
\end{eqnarray*}

Hence, to prove Theorem 2.1 (a), it suffices to show Theorem 2.2.

\textit{(b)} It can be proved analogously. To complete the proof of Theorem
2.1. (b), it suffices to express Theorem 2.3.

\bigskip

{\large \textbf{Theorem 2.2.}} \textit{The potential function }$V\in
L_{1}(0,\pi )$ \textit{of the problem (1.1),(1.7) satisfies}%
\begin{equation*}
V(x)=\lim_{n\rightarrow \infty }\left[ l_{n}^{j,1}\lambda _{n}-m^{2}\frac{%
l_{n}^{j,1}}{2\lambda _{n}}-\pi \right] \lambda _{n},
\end{equation*}%
\textit{for almost every} $x\in (0,\pi )$, \textit{with} $j_{1}=j_{n,1}(x).$

\bigskip

{\large \textbf{Proof. }}Lemma 1.3 yields%
\begin{equation*}
l_{n}^{j,1}-\frac{m^{2}l_{n}^{j,1}}{2\lambda _{n}^{2}}=\frac{\pi }{\lambda
_{n}}+\frac{\dint\limits_{x_{n}^{j,1}}^{x_{n}^{j+1,1}}V(t)dt}{\lambda _{n}}%
+O\left( \frac{1}{\lambda _{n}^{3}}\right) ,
\end{equation*}%
so that 
\begin{equation}
\left[ \left( l_{n}^{j,1}-\frac{m^{2}l_{n}^{j,1}}{2\lambda _{n}^{2}}\right)
\lambda _{n}-\pi \right] \lambda _{n}=\lambda
_{n}\dint\limits_{x_{n}^{j,1}}^{x_{n}^{j+1,1}}V(t)dt+O\left( \frac{1}{%
\lambda _{n}}\right) .  \tag{2.1}
\end{equation}

We may assume $x_{n}^{j,1}\neq x.$ By Theorem 1.1., if we take limit of both
sides of (2.1) as $n\rightarrow \infty $ for almost $x\in (0,\pi ),$ we get%
\begin{equation*}
V(x)=\lim_{n\rightarrow \infty }\left[ l_{n}^{j,1}\lambda _{n}-m^{2}\frac{%
l_{n}^{j,1}}{2\lambda _{n}}-\pi \right] \lambda _{n}.
\end{equation*}

{\large \textbf{Theorem 2.3.}} \textit{The potential function }$\widetilde{V}%
\in L_{1}(0,\pi )$ \textit{of the problem (1.1),(1.12) satisfies}%
\begin{equation*}
\widetilde{V}(x)=\lim_{n\rightarrow \infty }\left[ \widetilde{l}_{n}^{j,1}%
\widetilde{\lambda }_{n}\pm m^{2}\frac{(\widetilde{x}_{n}^{j,1}+\widetilde{x}%
_{n}^{j+1,1})}{2\widetilde{\lambda }_{n}}-\pi \right] \widetilde{\lambda }%
_{n}\pm m\sin 2\widetilde{\alpha },
\end{equation*}%
\textit{for almost every} $x\in (0,\pi )$, \textit{with} $j_{1}=j_{n,1}(x).$

\bigskip

{\large \textbf{Proof: }}It can be proved by using similar process to
Theorem 2.2.

\bigskip

{\large 3}\textbf{{\large . \textbf{M}ain Results }}

In this section, we solve a Lipschitz stability problem for Dirac operator.
Lipschitz stability is about a continuity between two metric spaces. So, we
have to first construct these spaces. To show continuity, we use a
homeomorphism between these spaces. Stability problems were studied by many
authors \cite{Cheng,Mrchenko,Mclghln}. To solve stability problem, we give a
main theorem which execute that the inverse nodal problem for Dirac system
is stable with Lipschitz stability. Here and later, we denote the space of
all admissible nodal sequences which converge to $V$ by $X=\left\{
X_{n}^{k,i}\right\} $ where $L_{n}^{k,i}=X_{n}^{k+1,i}-X_{n}^{k,i},i=1,2.$

\bigskip

\textbf{{\large \textbf{Definition} 3.1. }}\textit{Let} $%
\mathbb{N}
^{\prime }=%
\mathbb{N}
-\{1\}.$ \textit{We denote the space }$\Omega _{Dir}$ \textit{of all
potential functions of Dirac system and the space} $\Sigma _{Dir}$ \textit{%
of all admissible sequences by}

\begin{description}
\item[(i)] 
\begin{equation*}
\Omega _{Dir}=\{V\in L_{1}[0,\pi ]:V\text{ \textit{is the potential function
of the Dirac system}}\mathit{\},}
\end{equation*}%
\textit{and} $\Sigma _{Dir}=$\textit{The collection of the all double
sequences defined as}%
\begin{equation*}
X=\left\{ X_{n}^{k,1}:k=1,2,...,n;n\in 
\mathbb{N}
^{\prime },0<X_{n}^{1,1}<X_{n}^{2,1}<...<X_{n}^{n-1,1}<\pi \right\}
\end{equation*}%
\textit{for each} $n\in 
\mathbb{N}
.$

\item[(ii)] \textit{Let} $X\in \Sigma _{Dir}$ \textit{and define} $%
X=\{X_{n}^{k,1}\}$ \textit{where} $I_{n}^{k,1}=\left(
X_{n}^{k,1},X_{n}^{k+1,1}\right) .$ \textit{We say} $X$ \textit{is
quasinodal to some} $V\in L_{1}(0,\pi )$ \textit{if} $X$ \textit{is an
admissible sequence of nodes and satisfies} \textit{(I) and (II) below:}
\end{description}

\textit{(I)} $X$ has the following asymptotics uniformly for $k,$ as $%
n\rightarrow \infty $

\begin{equation*}
X_{n}^{k,1}=\frac{k\pi }{n-2}+O\left( \frac{1}{n}\right) ,k=1,2,...,n
\end{equation*}%
for the problem (1.1), (1.7). And the sequence%
\begin{equation*}
F_{n}=(n-2)\left\{ \dsum\limits_{k=1}^{n-1}\left[ n-2-\frac{m^{2}}{2(n-2)}%
\right] L_{n}^{k,1}-\pi \right\} ,
\end{equation*}%
converges to $V$ in $L_{1}.$

\textit{(II)} For the problem (1.1), (1.12), $X$ has below asymptotics
uniformly for $k,$ as $n\rightarrow \infty $

\begin{equation*}
X_{n}^{k,1}=\frac{k\pi }{n}+O\left( \frac{1}{n}\right) ,k=1,2,...,n
\end{equation*}%
and the sequence%
\begin{equation*}
F_{n}=n\left\{ \dsum\limits_{k=1}^{n-1}\left[ n\pm \frac{m^{2}}{2n}\right]
L_{n}^{k,1}\pm \frac{m^{2}}{n}X_{n}^{k,1}-\pi \right\} ,
\end{equation*}%
converges to $V$ in $L_{1}$. $X\in \Sigma _{Dir}$ is nodal if $X$ satisfies
one of the above asymptotic behaviours.

\bigskip

We denote $\Omega _{Dir}$ as a collection of all Dirac operators and the
space $\Sigma _{Dir}$ as a collection of all admissible double sequences of
nodes such that related functions are convergent in $L_{1}.$ A pseudometric $%
d_{\Sigma _{Dir}}$ on $\Sigma _{Dir}$ will be defined$.$ For convenience, we
will use the notation $X$ for the first component$.$ Essentially, $d_{\Sigma
_{Dir}}(X,\overline{X})$ is so close to%
\begin{equation*}
d_{0}(X,\overline{X})=\overline{\lim_{n\rightarrow \infty }}\pi \left[ n-2-%
\frac{m^{2}}{2(n-2)}\right] \dsum\limits_{k=1}^{n-1}\left\vert L_{n}^{k,1}-%
\overline{L}_{n}^{k,1}\right\vert ,
\end{equation*}%
where $n>\dfrac{m}{\sqrt{2}}+2$ and $L_{n}^{k,1}=X_{n}^{k+1,1}-X_{n}^{k,1}$, 
$\overline{L}_{n}^{k,1}=\overline{X}_{n}^{k+1,1}-\overline{X}_{n}^{k,1}.$

If we define $X\sim \overline{X}$ if and only if $d_{\Sigma _{Dir}}(X,%
\overline{X})=0,$ then $\sim $ is an equivalence relation on $\Sigma _{Dir}$
and $d_{\Sigma _{Dir}}$ would be a metric for the partition set $\Sigma
_{Dir}^{\ast }=\Sigma _{Dir}/\sim .$ Let $\Sigma _{Dir_{1}}\subset \Sigma
_{Dir}$ be the subspace of all asymptotically equivalent nodal sequences and
let $\Sigma _{Dir_{1}}^{\ast }=\Sigma _{Dir_{1}}/\sim .$ Let $\Phi $ be a
homeomorphism the maps $\Omega _{Dir}$ onto $\Sigma _{Dir_{1}}^{\ast }.$ We
will call $\Phi $ as a nodal map.

\textbf{{\large \textbf{Lemma} 3.1. }}\textit{Let} $X,\overline{X}\in \Sigma
_{Dir}.$

\bigskip

\textit{a)} \textit{If} $X$ \textit{belongs to case I, then}%
\begin{equation*}
L_{n}^{k,1}=\frac{\pi }{n-2}+O\left( \frac{1}{n}\right) ,k=1,2,...,n.
\end{equation*}

\textit{If} $X$ \textit{belongs to case II, then}%
\begin{equation*}
L_{n}^{k,1}=\frac{\pi }{n}+O\left( \frac{1}{n}\right) ,k=1,2,...,n.
\end{equation*}

\textit{b)} $\chi _{n,k}=\left\vert X_{n}^{k,1}-\overline{X}%
_{n}^{k,1}\right\vert =O\left( \frac{1}{n}\right) .$

\textit{c) For all }$x\in (0,\pi ),$ \textit{define} $J_{n,1}(x)=maks%
\{k:X_{n}^{k,1}\leq x\}$ \textit{so that} $k=J_{n,1}(x)$ \textit{if and only
if} $x\in \lbrack X_{n}^{J_{1},1},X_{n}^{J_{1}+1,1}].$ \textit{Then, for
sufficiently large} $n,$%
\begin{equation*}
\left\vert J_{n,1}(x)-\overline{J}_{n,1}(x)\right\vert \leq 1.
\end{equation*}

\textbf{{\large \textbf{Proof:} }}

\textit{a)} For case I, we get%
\begin{equation*}
L_{n}^{k,1}=X_{n}^{k+1,1}-X_{n}^{k,1}=\frac{\pi }{n-2}+O\left( \frac{1}{n}%
\right) ,
\end{equation*}%
by using the definition of nodal lengths. Similarly, for the case (II), we
obtain%
\begin{equation*}
L_{n}^{k,1}=X_{n}^{k+1,1}-X_{n}^{k,1}=\frac{\pi }{n}+O\left( \frac{1}{n}%
\right) .
\end{equation*}

\textit{b)} We only consider case I. The other case is similar. By using
asymptotic estimates, we get%
\begin{equation*}
\left\vert \chi _{n,k}\right\vert =\left\vert X_{n}^{k,1}-\overline{X}%
_{n}^{k,1}\right\vert \leq \left\vert X_{n}^{k,1}-\frac{k\pi }{n-2}%
\right\vert +\left\vert \frac{k\pi }{n-2}-\overline{X}_{n}^{k,1}\right\vert
=O\left( \frac{1}{n}\right) .
\end{equation*}

\textit{c) }Fix $x\in (0,\pi ).$ Let $J_{1}=J_{n,1}(x)$ and $\overline{J}%
_{1}=\overline{J}_{n,1}(x).$ Since%
\begin{equation*}
X_{n}^{J_{1},1}\leq x\leq X_{n}^{J_{1}+1,1}\Rightarrow \frac{J_{1}\pi }{n-2}%
+O\left( \frac{1}{n}\right) =X_{n}^{J_{1},1}\leq x\leq X_{n}^{J_{1}+1,1}=%
\frac{(J_{1}+1)\pi }{n-2}+O\left( \frac{1}{n}\right) ,
\end{equation*}%
and%
\begin{equation*}
\overline{X}_{n}^{\overline{J}_{1},1}\leq x\leq \overline{X}_{n}^{\overline{J%
}_{1}+1,1}\Rightarrow \frac{\overline{J}_{1}\pi }{n-2}+O\left( \frac{1}{n}%
\right) =\overline{X}_{n}^{\overline{J}_{1},1}\leq x\leq \overline{X}_{n}^{%
\overline{J}_{1}+1,1}=\frac{\overline{J}_{1}\pi }{n-2}+O\left( \frac{1}{n}%
\right) ,
\end{equation*}%
when $n$ is large enough, $\overline{J}_{1}+1\geq J_{1}$ and $J_{1}+1\geq 
\overline{J}_{1}.$ Hence, $-1\leq \overline{J}_{1}-J_{1}\leq 1,$ then $%
\left\vert \overline{J}_{1}-J_{1}\right\vert \leq 1.$

\bigskip

\textbf{{\large \textbf{Definition} 3.2. }}Suppose that $X,\overline{X}\in
\Sigma _{Dir}$ with $L_{n}^{k,1}$ and $\overline{L}_{n}^{k,1}$ are their
respective grid lengths. Let%
\begin{equation}
S_{n}\left( X,\overline{X}\right) =\pi \left[ n-2-\frac{m^{2}}{2(n-2)}\right]
\dsum\limits_{k=1}^{n-1}\left\vert L_{n}^{k,1}-\overline{L}%
_{n}^{k,1}\right\vert .  \tag{3.1}
\end{equation}%
Define%
\begin{equation*}
d_{0}\left( X,\overline{X}\right) =\overline{\lim\limits_{n\rightarrow
\infty }}S_{n}\left( X,\overline{X}\right) \text{ and }d_{\Sigma
_{Dir}}\left( X,\overline{X}\right) =\overline{\lim\limits_{n\rightarrow
\infty }}\frac{S_{n}\left( X,\overline{X}\right) }{1+S_{n}\left( X,\overline{%
X}\right) }.
\end{equation*}

We get this metric by evaluating $\left\Vert V-\overline{V}\right\Vert _{1}$
in Theorem 3.1. This definition was first made by \cite{Tsay}. Since the
function $f(x)=\dfrac{x}{1+x}$ is monotonic, we have 
\begin{equation*}
d_{\Sigma _{Dir}}\left( X,\overline{X}\right) =\frac{d_{0}\left( X,\overline{%
X}\right) }{1+d_{0}\left( X,\overline{X}\right) }\in \left[ 0,\pi \right] ,
\end{equation*}%
admitting that if $d_{0}\left( X,\overline{X}\right) =\infty ,$ then $%
d_{\Sigma _{Dir}}\left( X,\overline{X}\right) =1.$ Conversely 
\begin{equation*}
d_{0}\left( X,\overline{X}\right) =\frac{d_{\Sigma _{Dir}}\left( X,\overline{%
X}\right) }{1-d_{\Sigma _{Dir}}\left( X,\overline{X}\right) }.
\end{equation*}%
We can easily prove this equality by using Law's method \cite{Tsay,Cheng}.

\bigskip

\textbf{{\large \textbf{Lemma} 3.2. }}\textit{Let} $X,\overline{X}\in \Sigma
_{Dir}.$

\textit{a)} $d_{\Sigma _{Dir}}$ \textit{is a pseudometric on} $\Sigma
_{Dir}. $

\textit{b)} \textit{If} $X$\textit{\ and} $\overline{X}$ \textit{belong to
different cases,} $d_{\Sigma _{Dir}}\left( X,\overline{X}\right) =1.$

\bigskip

\textbf{\large Proof: } It can be proved similar way with in \cite{Tsay}.

\bigskip

Stability problems for Sturm-Liouville and Hill's operators were studied in 
\cite{Tsay,Cheng} respectively. Now, we prove the stability of the inverse
nodal problem for Dirac operator with Lipschitz stability. The below theorem
guarantees the Lipschitz stability of inverse nodal problem for Dirac
operator.

\bigskip

{\large \textbf{Theorem 3.1.}}\textit{\ The metric spaces }$\left( \Omega
_{Dir},\left\Vert .\right\Vert _{1}\right) $ and $\left( \Sigma
_{Dir_{1}}^{\ast }/\sim ,d_{\Sigma _{Dir}}\right) $ \textit{are homeomorphic
to each other. Here,} $\sim $ \textit{is the equivalence relation induced by}
$d_{\Sigma _{Dir}}.$ \textit{Furthermore} 
\begin{equation*}
\left\Vert V-\overline{V}\right\Vert _{1}=\frac{d_{\Sigma _{Dir}}\left( X,%
\overline{X}\right) }{1-d_{\Sigma _{Dir}}\left( X,\overline{X}\right) },
\end{equation*}%
\textit{where} $d_{\Sigma _{Dir}}\left( X,\overline{X}\right) <1$\textit{.}

\bigskip

{\large \textbf{Proof: }}By Lemma 3.2. we only need to consider when $X,%
\overline{X}\in \Sigma _{Dir}$ belong to same case. Without loss of
generality, let $X,\overline{X}$ belong to case I. In this case, we should
denote%
\begin{equation*}
\left\Vert V-\overline{V}\right\Vert _{1}=d_{0}\left( X,\overline{X}\right) .
\end{equation*}

According to the Theorem 2.1., $F_{n}$ and $\overline{F}_{n}$ converge to $V$
and $\overline{V},$ respectively. If we use the definition of norm in $L_{1}$
for the functions $V$ and $\overline{V},$ we have

\begin{equation*}
\left\Vert V-\overline{V}\right\Vert _{1}=\left( n-2\right) \left[ n-2-\frac{%
m^{2}}{2(n-2)}\right] \dint\limits_{0}^{\pi }\left\vert L_{n}^{J_{n,1}(x),1}-%
\overline{L}_{n}^{\overline{J}_{n,1}(x),1}\right\vert dx+o(1).
\end{equation*}

Hence by Fatou's Lemma,%
\begin{eqnarray}
\left\Vert V-\overline{V}\right\Vert _{1} &\leq &\left( n-2\right) \left[
n-2-\frac{m^{2}}{2(n-2)}\right] \dint\limits_{0}^{\pi }\left\vert
L_{n}^{J_{n,1}(x),1}-\overline{L}_{n}^{J_{n,1}(x),1}\right\vert dx  \notag \\
&&+\left( n-2\right) \left[ n-2-\frac{m^{2}}{2(n-2)}\right]
\dint\limits_{0}^{\pi }\left\vert \overline{L}_{n}^{J_{n,1}(x),1}-\overline{L%
}_{n}^{\overline{J_{n,1}}(x),1}\right\vert dx.  \TCItag{3.2}
\end{eqnarray}

Here, first and second terms can be written as

\begin{equation*}
\dint\limits_{0}^{\pi }\left\vert L_{n}^{J_{n,1}(x),1}-\overline{L}%
_{n}^{J_{n,1}(x),1}\right\vert dx=\frac{\pi }{n-2}\dsum\limits_{k=1}^{n-1}%
\left\vert L_{n}^{k,1}-\overline{L}_{n}^{k,1}\right\vert +o\left( \frac{1}{%
n^{2}}\right) ,
\end{equation*}%
and%
\begin{equation*}
\dint\limits_{0}^{\pi }\left\vert \overline{L}_{n}^{J_{n,1}(x),1}-\overline{L%
}_{n}^{\overline{J_{n,1}}(x),1}\right\vert dx=o\left( \frac{1}{n^{3}}\right)
.
\end{equation*}

If we consider last equalities in (3.2), we get%
\begin{eqnarray*}
\left\Vert V-\overline{V}\right\Vert _{1} &\leq &\left( n-2\right) \left[
n-2-\frac{m^{2}}{2(n-2)}\right] o\left( \frac{1}{n^{3}}\right) \\
&&+\left( n-2\right) \left[ n-2-\frac{m^{2}}{2(n-2)}\right] \left[ \frac{\pi 
}{n-2}\dsum\limits_{k=1}^{n-1}\left\vert L_{n}^{k,1}-\overline{L}%
_{n}^{k,1}\right\vert +o\left( \frac{1}{n^{2}}\right) \right] ,
\end{eqnarray*}%
and%
\begin{equation}
\left\Vert V-\overline{V}\right\Vert _{1}\leq \pi \left[ n-2-\frac{m^{2}}{%
2(n-2)}\right] \dsum\limits_{k=1}^{n-1}\left\vert L_{n}^{k,1}-\overline{L}%
_{n}^{k,1}\right\vert +o(1).  \tag{3.3}
\end{equation}

Similarly, we can get

\begin{equation}
\left\Vert V-\overline{V}\right\Vert _{1}\geq \pi \left[ n-2-\frac{m^{2}}{%
2(n-2)}\right] \dsum\limits_{k=1}^{n-1}\left\vert L_{n}^{k,1}-\overline{L}%
_{n}^{k,1}\right\vert +o(1).  \tag{3.4}
\end{equation}

If we consider (3.3) and (3.4) together, we obtain%
\begin{equation*}
\left\Vert V-\overline{V}\right\Vert _{1}=\pi \left[ n-2-\frac{m^{2}}{2(n-2)}%
\right] \dsum\limits_{k=1}^{n-1}\left\vert L_{n}^{k,1}-\overline{L}%
_{n}^{k,1}\right\vert .
\end{equation*}

The proof is complete after taking the limit as $n\rightarrow \infty $.

\bigskip

\bigskip

\textbf{\large 4. Conclusion}

In this study, Lipschtiz stability of inverse nodal problem is proved for
Dirac operator by using zeros of the first component function of two
dimensional vector eigenfunction. Especially, two metric spaces were defined
and it was shown that they were homeomorphic to each other. These results
are new and can be generalized.

\bigskip

\textbf{Acknowledgements}

The authors are grateful to the reviewers for their valuable comments that
greatly improved the paper.

\bigskip 

\textbf{Conflict of interest}

On behalf of all authors, the corresponding author states that there is no
conflict of interest.

\end{document}